\documentclass[12pt,reqno]{amsart}
\usepackage[margin=1in]{geometry}
\usepackage{amsmath,amssymb,amsthm,graphicx,amsxtra, setspace}
\usepackage[utf8]{inputenc}
\usepackage{mathrsfs}
\usepackage{hyperref}
\usepackage{xcolor}
\usepackage{upgreek}
\usepackage{mathtools}
\usepackage[numbers]{natbib}
\allowdisplaybreaks
\newtheorem{thm}{Theorem}[section]

\newtheorem{lem}{Lemma}[section]
\newtheorem{rem}{Remark}[section]

\newtheorem{Def}{Definition}[section]

\newtheorem{Ass}{Assumption}[section]
\usepackage{latexsym}
\usepackage{hyperref}
\usepackage{graphicx}
\usepackage{epstopdf}

\allowdisplaybreaks

\let\originalleft\left
\let\originalright\right
\renewcommand{\left}{\mathopen{}\mathclose\bgroup\originalleft}
\renewcommand{\right}{\aftergroup\egroup\originalright}

\newcommand{\Addresses}{{
		\footnote{

			\noindent	 \textsuperscript{1,2} Department of Applied Mathematics and Scientific Computing, Indian Institute of Technology Roorkee, Roorkee, 247667, India.	
			
			\noindent  \textit{e-mail\textsuperscript{1}:} \texttt{g\_gupta@as.iitr.ac.in}
			
			\noindent  \textit{e-mail\textsuperscript{2}:} \texttt{jay.dabas@gmail.com.}

			\noindent \textsuperscript{*}Corresponding author.
			
			\textit{Key words:} Non-autonomous differential equations, impulsive systems, Approximate Controllability, Krasnoselskii's fixed point theorem.
			
			Mathematics Subject Classification (2010): 35A12; 37L05; 93C27; 93B05; 93C10. 
			
}}}

\begin{document}
	\title[]{The existence and controllability of nonautonomous system influenced by impulses on both state and control.\Addresses}
	\author [Garima Gupta AND Jaydev Dabas]{Garima Gupta\textsuperscript{1} AND  Jaydev Dabas\textsuperscript{2*}}
	\maketitle

	\begin{abstract}
		This paper examines impulsive controls related to nonautonomous impulsive integro-differential equations in Hilbert space, highlighting their significance. We establish the existence of the mild solution by using fixed point approach and present conditions for approximate controllability using impulsive resolvent operators and the adjoint problem, supported by an illustrative example.
	\end{abstract}
	\section{\textbf{Introduction}}\label{intro}\setcounter{equation}{0}
	\vspace{0.5cm}
	
		Controllability is a fundamental concept in mathematical control theory, crucial for addressing numerous control issues, including stabilizing unstable systems via feedback control\cite{barbu2018controllability}, ensuring the irreducibility of transition semigroups\cite{da1996ergodicity}, and solving optimal control problems\cite{barbu1993analysis}. A system is considered controllable if every state within the system can be guided to a desired outcome within a specified time frame using appropriate control inputs. Various forms of controllability have been developed, such as exact, null, approximate, interior, boundary, and finite-approximate controllability \cite{zuazua2006controllability}. In the context of infinite-dimensional control systems, approximate controllability plays a significant role due to its wide range of applications, see\cite{triggiani1977note}.
		
	    The theory of impulsive differential equations offers a significantly broader scope compared to the theory of differential equations without impulse effects. Even a basic impulsive differential equation can give rise to unique behaviors, such as rhythmic oscillations, solution merging, and the noncontinuability of solutions. Several articles has been published considering impusive systems. Recently in 2024, Mahmudov studied the approximate controllability of the following impulsive system in a Hilbert space $\mathbb{H}$ \cite{mahmudov2024study}
			\begin{align}\label{lip}
		\begin{split}
			x^{\prime}(t)&=\mathrm{A} x(t)+\mathrm{B} u(t), \quad t \in J=[0, b] \backslash\left\{t_{1}, \ldots, t_{m}\right\}, \\ 
			\Delta x\left(t_{k}\right)&=\mathrm{D}_{k} x\left(t_{k}\right)+\mathrm{E}_{k} v_{k}, \quad k=1, \ldots, m, \\ 
			x(0)&=x_{0}. 
		\end{split}			
	\end{align}
	Here, $x(.)\in\mathbb{H}$ with norm $\|x\|=\sqrt{\langle x,x\rangle}$ and control function $u(.)\in L^{2}\left([0,b],\mathbb{U}\right)$, where $\mathbb{U}$ is another Hilbert space, and $v_{k}\in \mathbb{U} \text{ for } k=1,2,\dots m$.
	The study of these type of impulsive systems has been done in finite dimensional cases before (for reference one can read \cite{leela1993controllability,george2000note,zhao2010controllability,han2013note,s2020controllability}).
	The work by Mahmudov described by equation\ref{lip} is the first work with these type of impulse structure in infinite dimensional space. In that study, Mahmudov developed the solution by employing semigroups and impulsive operators, and he derived both the necessary and sufficient conditions for the approximate controllability of linear impulsive evolution equations using the concept of the impulsive resolvent operator.\newline
	From a mathematical perspective, the impulsive effect  $\Delta x\left(t_{k+1}\right)=\mathrm{D}_{k+1}x\left(t_{k+1}\right)+\mathrm{E}_{k+1}v_{k+1}$
	signifies that at each impulsive time $t_{k+1}$ the state $x\left(t_{k+1}\right)$  undergoes an immediate change due to a linear transformation $\mathrm{D}_{k+1}$ combined with an additional input $\mathrm{E}_{k+1}v_{k+1}$. This method is especially useful for capturing scenarios where the state experiences abrupt shifts caused by external forces or internal dynamics.
	
	In classical impulsive systems, the impulsive effects are generally characterized by abrupt changes in the system's state at specified time instances, typically represented by a jump condition of the form:
	\begin{align*}
		\Delta x\left(x(t_{k})\right) = I_{k}\left(x\left(t_{k}^{-}\right)\right),
	\end{align*}
	where $I_{k}$ is the impulse function, and $x\left(t_{k}^{-}\right)$ indicates the system's state just before the impulse occurs at time $t_{k}$. 
	
	These conventional impulsive mechanisms are suitable for describing systems influenced by regular or expected disturbances. However, they might not accurately depict the complexities of systems that are subject to both foreseeable and unforeseeable impulses. The impulse function $I_{k}$ usually represents a predetermined modification to the state, which can be overly restrictive when dealing with scenarios involving intricate or unpredictable changes.
	
	As a result, while the traditional impulse model works $\Delta x\left(x(t_{k})\right)= I_{k}\left(x\left(t_{k}^{-}\right)\right)$ works well for simpler or more predictable systems, the alternative approach $\Delta x\left(t_{k+1}\right)=\mathrm{D}_{k+1}x\left(t_{k+1}\right)+\mathrm{E}_{k+1}v_{k+1}$  offers a more sophisticated and flexible framework for representing frequent and intricate abrupt shifts in a system’s state.
	
	The next work in this direction is being done by Javed A Asadzade and Mahmudov. They have studied the existence, optimal control of impulsive sctochastic evolution systems\cite{asadzade2024solvability}. They have also studied the approximate controllability of semilinear systems\cite{asadzade2024approximate}.
	
	The study of approximate controllability in nonautonomous systems is vital for managing time-dependent dynamics in real-world applications. It enables systems to follow desired behaviors despite uncertainties, using adaptive and flexible controls. This is crucial in fields like robotics (following time-varying paths), climate science (modeling seasonal changes), economics (managing market fluctuations), and biomedical engineering (targeted therapies). For the detailed study of nonautonomous systems we refer a  book by Kloeden et al.\cite{kloeden2011nonautonomous}. The articles by Arora et al.\cite{arora2022existence} and Ravikumar et.al are really appreciable for the study of approximate controllability in nonautonomous systems.
	
	From the above literature it is clear that the study of approximate controllability of semilinear deterministic system with the specified impulse as in system \eqref{lip} has not been studied yet and motivate us to consider the nonautonomous integro impulsive systems within a separable Hilbert space as follows:
	\begin{align}\label{1}
		\begin{split}
			x^{\prime}(t)&=\mathrm{A}(t) x(t)+\mathrm{B} u(t)+f\left(t,x(t)\right) +\int_{0}^{t}q(t-s) \xi(s,x(s)) ds, \quad t \in J=[0, b] \backslash\left\{t_{1}, \ldots, t_{m}\right\}, \\ 
			\Delta x\left(t_{k}\right)&=\mathrm{D}_{k} x\left(t_{k}\right)+\mathrm{E}_{k} v_{k}, \quad k=1, \ldots, m, \\ 
			x(0)&=x_{0}, 
		\end{split}			
	\end{align}
	
	where $\{\mathrm{A}(t) : t\in J\}$ is a family of linear operators (not necessarily bounded) on $\mathbb{H}$.
	Control $u:J\rightarrow\mathbb{U}$, where $\mathbb{U}$ is Hilbert space identified with its own dual is given in $L^{2}([0, b], \mathbb{U})$, $v_{k} \in \mathbb{U}, k=1, \ldots, m$. $\mathrm{B}:\mathbb{U}\rightarrow\mathbb{H}$, $\mathrm{D}_{k}:\mathbb{H}\rightarrow\mathbb{H}$, $\mathrm{E}_{k}:\mathbb{U}\rightarrow\mathbb{H}$ are bounded linear operators and $\|\mathrm{B}\|_\mathcal{L}=M_{B}$. The functions
	$f, \xi:J\times \mathbb{H}\rightarrow\mathbb{H}$ are satisfying some suitable assumptions. $q: [0,b]\rightarrow \mathbb{H}$ is continuous and $q\in L^{1}\left([0,b],\mathbb{R}^{+}\right)$.
	
	At the points of discontinuance $t_{k}$ $(\text{where}\quad k=1,\dots, m \quad\text{and}\quad 0=t_0<t_1<t_2<\dots<t_m<t_{n+1}=b )$, 
	the state variable's abrupt change is determined by $\Delta x\left(t_{k}\right)=x(t_{k}^{+})-x({t_{k}^-})$, with $x(t_{k}^{\pm})=\lim_{h\to0^{\pm}}x(t_{k}+h)$ and the supposition that $x({t_{k}^-})=x({t_{k}})$.
	
   $\prod_{j=1}^{k}A_{j}$ denotes the operator composition $A_{1}A_{2}\dots A_{k}$. For $j=k+1$ to $k$, $\prod_{j=k}^{k+1}A_{j}=1$. In the same way, $\prod_{j=k}^{1}A_{j}$ represents the composition $A_{k}A_{k-1}\dots A_{1}$ and $\prod_{j=k+1}^{k}A_{j}=1$.
	
	\section{\textbf{Preliminaries and Assumptions}}\label{pre}\setcounter{equation}{0}
	This section contains some essential definitions and specified assumptions which are required to derive the sufficient conditions for ensuring the approximate controllability of system \eqref{1}.
	
	Let us introduce
	\begin{align*}
		\mathcal{PC}(J;\mathbb{H})&:=\{\psi:J\rightarrow\mathbb{H}: \psi(\cdot)\ \mbox{is piecewise continuous with jump}\\ \quad & \mbox{discontinuity at}\ t_{k}\ \mbox{satisfying}\ x(t_{k}^{-})=x(t_{k})\}.
	\end{align*}
	For $x\in \mathcal{PC}(J;\mathbb{H})$, we define $\|x\|_\mathcal{PC}=\sup_{t\in J}\|x(t)\|$.
	
	\subsection{Evolution family} An evolution family is an essential concept in the study of non-autonomous systems, particularly when dealing with time-dependent differential equations. Here is a formal definition:
	
\begin{Def}[\cite{pazy2012semigroups}]
		Let $\mathbb{X}$ be a Banach space, and let $J=[0,b]$, be an interval of the real line. An evolution family $\{ \mathrm{U}(t,s) \}_{(t,s) \in J \times J, \, t \geq s}$ is a two-parameter family of bounded linear operators on $\mathbb{X}$ with the following properties:
	
	\begin{enumerate}
		\item Initial Condition:
		\[
		\mathrm{U}(s,s) = \mathrm{I} \quad \text{for all } s \in J,
		\]
		where $\mathrm{I}$ is the identity operator on $\mathbb{X}$.
		
		\item Semigroup Property (also called the cocycle condition):
		\[
		\mathrm{U}(t,s) = \mathrm{U}(t,r) \mathrm{U}(r,s) \quad \text{for all } s \leq r \leq t \text{ in } J.
		\]
		
		\item Strong Continuity: The mapping $(t,s) \mapsto \mathrm{U}(t,s)x$ is continuous for each fixed $x \in \mathbb{X}$.
	\end{enumerate}
\end{Def}
	
	To construct an evolution family, let us impose the following assumptions on the family of linear operators $\{\mathrm{A}(t) : t\in J\}$ (see, chapter 5,\cite{pazy2012semigroups}).	

\begin{Ass}\label{2.1}
	\begin{enumerate}
		\item[(R1)] The linear operator $\mathrm{A}(t)$ is closed for each $t\in J$ and the domain $\mathcal{D}(\mathrm{A}(t))=\mathcal{D}(\mathrm{A})$ is dense in $\mathbb{X}$ and independent of $t$.
		\item[(R2)] The resolvent operator $\mathrm{R}(\lambda, \mathrm{A}(t))$ for $t\in J$ exists for all $\lambda$ with $Re \lambda \le0$ and there exists $K>0$ such that 
		\begin{align*}
			\left\|\mathrm{R}(\lambda, \mathrm{A}(t))\right\|_{\mathcal{L}(\mathbb{X})}\le \frac{K}{|\lambda|+1}
		\end{align*}
	\item[(R3)] There exist constants $N>0$ and $0<\delta \le 1$ such that 
	\begin{align*}
		\left\|(\mathrm{A}(t)-\mathrm{A}(s))\mathrm{A}^{-1}(\tau)\right\|_{\mathcal{L}(\mathbb{X})}\le N |t-s|^{\delta}, \text{for all } t,s,\tau \in J.
	\end{align*}
\item[(R4)] The operator $\mathrm{R}(\lambda, \mathrm{A}(t)), t\in J$ is compact for some $\lambda \in \rho(\mathrm{A}(t))$, where $\rho(\mathrm{A}(t))$ is the resolvent set of $\mathrm{A}(t)$.
	\end{enumerate}
\end{Ass}
	
	\begin{lem} [Theorem 6.1, Chapter 5, \cite{pazy2012semigroups}]\label{lem 2.1}
		Suppose that (R1)-(R3) hold true. Then there esists a unique evolution family $\mathrm{U}(t,s)$ on $0\le s\le t\le b$ satisfying the following:
		\begin{enumerate}
			\item For $0\le s\le t\le b$, we have $\left\|\mathrm{U}(t,s)\right\|_{\mathcal{L}(\mathbb{X})}\le M$.
			\item The operator $\mathrm{U}(t,s): \mathbb{X}\mapsto \mathcal{D}(\mathrm{A})$ for $0\le s\le t\le b$ and the mapping $t \mapsto \mathrm{U}(t,s)$ is strongly differentiable in $\mathbb{X}$. The derivative $\frac{\partial}{\partial t}\mathrm{U}(t,s)\in \mathcal{L}(\mathbb{X})$ and it is strongly continuous on $0\le s\le t\le b$. Moreover, 
			\begin{align*}
				\frac{\partial}{\partial t}\mathrm{U}(t,s)-\mathrm{A}(t)\mathrm{U}(t,s)=0, \text{ for } 0\le s\le t\le b,
			\end{align*}
		\begin{align*}
			\left\|\frac{\partial}{\partial t}\mathrm{U}(t,s)\right\|_{\mathcal{L}(\mathbb{X})}=\left\|\mathrm{A}(t)\mathrm{U}(t,s)\right\|_{\mathcal{L}(\mathbb{X})}\le \frac{M}{t-s},
		\end{align*}
	and
	\begin{align*}
		\left\|\mathrm{A}(t)\mathrm{U}(t,s)\mathrm{A}(s)^{-1}\right\|_{\mathcal{L}(\mathbb{X})}\le M, \text{ for }  0\le s\le t\le b.
	\end{align*}
    \item For each $t\in J$ and every $v\in \mathcal{D}(\mathrm{A})$, $\mathrm{U}(t,s)v$ is differentiable with respect to $s$ on $0\le s\le t\le b$ and
    \begin{align*}
    	\frac{\partial}{\partial t}\mathrm{U}(t,s)v=-\mathrm{U}(t,s)\mathrm{A}(s)v.
    \end{align*}
		\end{enumerate}
	\end{lem}
	
	\begin{lem}[Proposition 2.1,\cite{fitzgibbon1978semilinear}]\label{lem2.2}
		Supose \{$\mathrm{A}(t): t\in J$\} satisfies the assumptions (R1)-(R4). Let \{$\mathrm{U}(t,s) :0\le s\le t\le b$\} be the linear evolution family generated by \{$\mathrm{A}(t): t\in J$\}, then \{$\mathrm{U}(t,s) :0\le s\le t\le b$\} is a compact operator, whenever $t-s >0$. 
	\end{lem}

	\begin{Def}
		A mild solution $x:J\rightarrow \mathbb{H}$ of the system \eqref{1} satisfying $x(0)=x_{0} $ and $\Delta x\left(t_{k}\right)=\mathrm{D}_{k} x\left(t_{k}\right)+\mathrm{E}_{k} v_{k}$, $k=1,\ldots,m$ on the intervals $t_{k-1}<t\le t_{k}$ is continuous, which is given by
		\begin{align}\label{2}
			x(t)=	\begin{cases}
				\mathrm{U}(t,0) x(0)+\int_{0}^{t} \mathrm{U}(t,s) [\mathrm{B} u(s)+f\left(s,x(s)\right)+\int_{0}^{s}q(s-\tau)\xi (\tau,x(\tau)) d\tau] d s, \quad 0 \leq t \leq t_{1} \\
				\mathrm{U}\left(t,t_{k}\right) x\left(t_{k}^{+}\right)+\int_{t_{k}}^{t} \mathrm{U}(t,s) [\mathrm{B} u(s)+f\left(s,x(s)\right)+\int_{0}^{s}q(s-\tau)\xi (\tau,x(\tau)) d\tau] d s,\\\quad t_{k}<t \leq t_{k+1},  k=1, \ldots, m,
			\end{cases}
		\end{align}
		with
		\begin{align*}
			\begin{split}
				x\left(t_{k}^{+}\right)  =&\prod_{j=k}^{1}\left(\mathrm{I}+\mathrm{D}_{j}\right) \mathrm{U}\left(t_{j},t_{j-1}\right) x_{0} +\sum_{i=1}^{k} \prod_{j=k}^{i+1}\left(\mathrm{I}+\mathrm{D}_{j}\right) \mathrm{U}\left(t_{j},t_{j-1}\right)\left(\mathrm{I}+\mathrm{D}_{i}\right) \\
				&\qquad \times\left(\int_{t_{i-1}}^{t_{i}} \mathrm{U}\left(t_{i},s\right) [\mathrm{B}u(s)+f\left(s,x(s)\right)+\int_{0}^{s}q(s-\tau)\xi (\tau,x(\tau)) d\tau] d s\right)\\& +\sum_{i=2}^{k} \prod_{j=k}^{i}\left(\mathrm{I}+\mathrm{D}_{j}\right)\mathrm{U}\left(t_{j},t_{j-1}\right) \mathrm{E}_{i-1} v_{i-1}+\mathrm{E}_{k} v_{k}. 
			\end{split}
		\end{align*}
	\end{Def}
	\begin{Def}\cite{mahmudov2003approximate}
		The system \eqref{1} is considered approximately controllable on the interval $J$ if the closure of the reachable set equals the entire space $\mathbb{H}$. The reachable set is defined by 
		\begin{align*}
			\mathfrak{R}_{t} = \{x \in \mathbb{H} \; | \; x = x(t,0,u), \; u(\cdot) \in \mathrm{L}^{2}(J; \mathbb{U})\}.
		\end{align*}
	\end{Def}
	
	\begin{lem}[Theorem 1, \cite{burton1998fixed}](Krasnoselskii's Fixed Point Theorem)\label{fixed point thm}
		Let $\mathcal{E}$ be a closed, bounded and convex subset of a Banach space $\mathbb{X}$ and let $\mathcal{G}_{1}$ and $\mathcal{G}_{2}$ be two mappings of $\mathcal{E}$ into $\mathbb{X}$ such that $\mathcal{G}_{1}(w)+\mathcal{G}_{2}(x)\in\mathcal{E}$, whenever $w,x\in\mathcal{E}$. If $\mathcal{G}_{1}$ is continuous and $\mathcal{G}_{1}(\mathcal{E})$ is relatively compact subset of $\mathcal{E}$ . Also $\mathcal{G}_{2}$ is a contraction map. Then there exists $z\in \mathcal{E}$ such that $z=\mathcal{G}_{1}(z)+\mathcal{G}_{2}(z)$.
	\end{lem}
	
		\subsection{Linear non-autonomous system}
	The linear nonautonomous impulsive system corresponding to system\eqref{1} in $\mathbb{H}$ is given by:
	\begin{align}\label{4}
		\begin{split}
			x^{\prime}(t)&=\mathrm{A}(t) x(t)+\mathrm{B} u(t), \quad t \in J=[0, b] \backslash\left\{t_{1}, \ldots, t_{m}\right\}, \\ 
			\Delta x\left(t_{k}\right)&=\mathrm{D}_{k} x\left(t_{k}\right)+\mathrm{E}_{k} v_{k}, \quad k=1, \ldots, m, \\ 
			x(0)&=x_{0}. 
		\end{split}			
	\end{align}
The mild solution of the above linear system is given by the following expression
\begin{align}\label{eqn2.3}
	x(t)=	\begin{cases}
		\mathrm{U}(t,0) x(0)+\int_{0}^{t} \mathrm{U}(t,s) \mathrm{B} u(s) d s, \quad 0 \leq t \leq t_{1} \\
		\mathrm{U}\left(t,t_{k}\right) x\left(t_{k}^{+}\right)+\int_{t_{k}}^{t} \mathrm{U}(t,s) \mathrm{B} u(s) d s,\\\quad t_{k}<t \leq t_{k+1},  k=1, \ldots, m,
	\end{cases}
\end{align}
with
\begin{align*}
	\begin{split}
		x\left(t_{k}^{+}\right)  =&\prod_{j=k}^{1}\left(\mathrm{I}+\mathrm{D}_{j}\right) \mathrm{U}\left(t_{j},t_{j-1}\right) x_{0} +\sum_{i=1}^{k} \prod_{j=k}^{i+1}\left(\mathrm{I}+\mathrm{D}_{j}\right) \mathrm{U}\left(t_{j},t_{j-1}\right)\left(\mathrm{I}+\mathrm{D}_{i}\right)\int_{t_{i-1}}^{t_{i}} \mathrm{U}\left(t_{i},s\right) \mathrm{B}u(s)\\& +\sum_{i=2}^{k} \prod_{j=k}^{i}\left(\mathrm{I}+\mathrm{D}_{j}\right)\mathrm{U}\left(t_{j},t_{j-1}\right) \mathrm{E}_{i-1} v_{i-1}+\mathrm{E}_{k} v_{k}. 
	\end{split}
\end{align*}
	To demonstrate the approximate controllability of the linear system mentioned above, we introduce a bounded linear operator $\mathrm{M}: L^{2}(J, \mathbb{U}) \times \mathbb{U}^{m} \rightarrow \mathbb{H}$ as follows:
	
	\begin{align*}
		& \mathrm{M}\left(u(\cdot),\left\{v_{k}\right\}_{k=1}^{m}\right) \\
		& =\mathrm{U}\left(b,t_{m}\right) \sum_{i=1}^{m} \prod_{j=m}^{i+1}\left(\mathrm{I}+\mathrm{D}_{j}\right) \mathrm{U}\left(t_{j},t_{j-1}\right)\left(\mathrm{I}+\mathrm{D}_{i}\right) \int_{t_{i-1}}^{t_{i}} \mathrm{U}\left(t_{i},s\right) \mathrm{B} u(s) d s \\
		&\qquad +\int_{t_{m}}^{b} \mathrm{T}(b-s) \mathrm{B} u(s) d s +\mathrm{U}\left(b,t_{m}\right) \sum_{i=2}^{m} \prod_{j=m}^{i}\left(\mathrm{I}+\mathrm{D}_{j}\right) \mathrm{U}\left(t_{j},t_{j-1}\right) \mathrm{E}_{i-1} v_{i-1}+\mathrm{U}\left(b,t_{m}\right) \mathrm{E}_{m} v_{m}.
	\end{align*}
\begin{rem}
	We can verify Lemma 7 and Lemma 9  in \cite{mahmudov2024study} for the linear system \ref{4} in similar way.
\end{rem}
	The operator $M^{*}$ is the adjoint of $M$ and has the following form (it can be verified in the similar way as in Lemma 9,\cite{mahmudov2024study})
	\begin{align*}
		\mathrm{M}^{*} \varphi & =\left(\mathrm{B}^{*} \psi(\cdot),\left\{\mathrm{D}_{k}^{*} \psi\left(t_{k}^{+}\right)\right\}_{k=1}^{m}\right), \\
		\mathrm{B}^{*} \psi(t) & = \begin{cases}\mathrm{B}^{*} \mathrm{U}^{*}(b,t) \varphi,\quad t_{m}<t \leq b, \\
			\mathrm{B}^{*} \mathrm{U}^{*}\left(t_{k},t\right)\left(\mathrm{I}+\mathrm{D}_{k}^{*}\right) \prod_{i=k+1}^{m} \mathrm{U}^{*}\left(t_{i},t_{i-1}\right)  \left(\mathrm{I}+\mathrm{D}_{i}^{*}\right) \mathrm{U}^{*}\left(b,t_{m}\right) \varphi, \quad t_{k-1}<t \leq t_{k},\end{cases} \\
		\mathrm{E}_{k}^{*} \psi\left(t_{k}^{+}\right) & = \begin{cases}\mathrm{E}_{m}^{*} \mathrm{U}^{*}\left(b,t_{m}\right) \varphi, \quad k=m, \\
			\mathrm{E}_{k}^{*} \prod_{i=k+1}^{m} \mathrm{U}^{*}\left(t_{i},t_{i-1}\right) \left(\mathrm{I}+\mathrm{D}_{i}^{*}\right) \mathrm{U}^{*}\left(b,t_{m}\right) \varphi, \quad k=m-1, \ldots, 1,\end{cases}
	\end{align*}
	where the operators $\mathrm{U}^{*}$, $\mathrm{B}^{*}$, $\mathrm{D}_{k}^{*}$, $\mathrm{E}_{m}^{*}$ are the adjoint operators of $\mathrm{U}$, $\mathrm{B}$, $\mathrm{D}_{k}$ and $\mathrm{E}_{k}$ respectively and $\psi(.)$ is the solution of the adjoint problem associated with system \ref{4}.
	The operator $\mathrm{M}\mathrm{M}^{*}:\mathbb{H}\rightarrow\mathbb{H}$ has the following form:
	\begin{align*}
		\mathrm{M}\mathrm{M}^{*}=\Theta_{0}^{t_{m}}+\Gamma_{t_{m}}^{b}+\widetilde{\Theta}_{0}^{t_{m}}+\widetilde{\Gamma}_{t_{m}}^{b},
	\end{align*}

	where $\Gamma_{t_{m}}^{b}, \widetilde{\Gamma}_{t_{m}}^{b}, \Theta_{0}^{t_{m}}, \widetilde{\Theta}_{0}^{t_{m}}: \mathbb{H} \rightarrow \mathbb{H}$ are non-negative operators and defined as follows:
	\begin{align*}
		\Gamma_{t_{m}}^{b}:= & \int_{t_{m}}^{b} \mathrm{U}(b,s) \mathrm{B} \mathrm{B}^{*} \mathrm{U}(b,s) d s, \quad \widetilde{\Gamma}_{t_{m}}^{b}:=\mathrm{U}\left(b,t_{m}\right) \mathrm{E}_{m} \mathrm{E}_{m}^{*} \mathrm{U}^{*}\left(b,t_{m}\right), \\
	\end{align*}
	
	\begin{align*}
		\Theta_{0}^{t_{m}}:= & \mathrm{U}\left(b,t_{m}\right)\\ &\times \sum_{i=1}^{m} \prod_{j=m}^{i+1}\left(\mathrm{I}+\mathrm{D}_{j}\right) \mathrm{U}\left(t_{j},t_{j-1}\right)\left(\mathrm{I}+\mathrm{D}_{i}\right) \int_{t_{i-1}}^{t_{i}} \mathrm{U}\left(t_{i},s\right) \mathrm{B}\mathrm{B}^{*} \mathrm{U}^{*}\left(t_{k},s\right) d s \\
		& \times\left(\mathrm{I}+\mathrm{D}_{i}^{*}\right) \prod_{k=i+1}^{m} \mathrm{U}^{*}\left(t_{k},t_{k-1}\right)\left(\mathrm{I}+\mathrm{D}_{k}^{*}\right) \mathrm{U}^{*}\left(b,t_{m}\right), \\
		\widetilde{\Theta}_{0}^{t_{m}}:= & \mathrm{U}\left(b,t_{m}\right) \sum_{i=2}^{m} \prod_{j=m}^{i}\left(I+\mathrm{D}_{j}\right) \mathrm{U}\left(t_{j},t_{j-1}\right) \mathrm{E}_{i-1} \mathrm{E}_{i-1}^{*} \\
		& \times \prod_{k=i}^{m} \mathrm{U}^{*}\left(t_{k},t_{k-1}\right)\left(I+\mathrm{D}_{k}^{*}\right) \mathrm{U}^{*}\left(b,t_{m}\right) .
	\end{align*}

    \begin{rem}
    	The linear system \ref{4} is said to be approximately controllable on $[0,b]$ if $\overline{\operatorname{Im} \mathrm{M}}=\mathbb{H}$
    \end{rem}
Now we will prove the approximate controllability of linear non-autonomous system\ref{4}.
 \begin{thm}
 	For the system \ref{4}, the following statements are equivalent:
 	\begin{enumerate}
 		\item[(a)] System \ref{4} is approximately controllable on $[0, b]$.
 		\item[(b)] $\mathrm{M}^{*} \varphi=0$ implies that $\varphi=0$.
 		\item[(c)]  $\Theta_{0}^{t_{m}}+\Gamma_{t_{m}}^{b}+\widetilde{\Theta}_{0}^{t_{m}}+\widetilde{\Gamma}_{t_{m}}^{b}$ is strictly positive.
 		\item[(d)]  $\lambda\left(\lambda \mathrm{I}+\Theta_{0}^{t_{m}}+\Gamma_{t_{m}}^{b}+\widetilde{\Theta}_{0}^{t_{m}}+\widetilde{\Gamma}_{t_{m}}^{b}\right)^{-1}$ converges to zero operator as $\lambda \rightarrow 0^{+}$ in strong operator topology.
 		\item[(e)] $\lambda\left(\lambda \mathrm{I}+\Theta_{0}^{t_{m}}+\Gamma_{t_{m}}^{b}+\widetilde{\Theta}_{0}^{t_{m}}+\widetilde{\Gamma}_{t_{m}}^{b}\right)^{-1}$ converges to zero operator as $\lambda \rightarrow 0^{+}$ in weak operator topology.
 	\end{enumerate}
 	 
 \end{thm}

Proof. The proof of the equivalence $(a) \Longleftrightarrow(b)$ is standard. Approximately controllability of system (1) on $[0, b]$ is equivalent to $\operatorname{Im} \mathrm{M}$ is dense in $\mathbb{H}$. That means, the kernel of $\mathrm{M}^{*}$ is trivial in $\mathbb{H}$. Equivalently,

\begin{align*}
	\mathrm{M}^{*} \varphi=\left(\mathrm{B}^{*} \psi(\cdot),\left\{\mathrm{E}_{k}^{*} \psi\left(t_{k}^{+}\right)\right\}_{k=1}^{m}\right)=0,
\end{align*}

implies that $\varphi=0$. For the equivalence $(a) \Longleftrightarrow(c)$ is clear from \cite{mahmudov2024study}. The equivalence $(d) \Longleftrightarrow(e)$ is a consequence of positivity of
\begin{align*}
	\lambda\left(\lambda \mathrm{I}+\Theta_{0}^{t_{m}}+\Gamma_{t_{m}}^{b}+\widetilde{\Theta}_{0}^{t_{m}}+\widetilde{\Gamma}_{t_{m}}^{b}\right)^{-1}.
\end{align*}
We prove only $(a) \Longleftrightarrow(d)$. To do so, consider the functional

\begin{align*}
	J_{\lambda}(\varphi)=\frac{1}{2}\left\|\mathrm{M}^{*} \varphi\right\|^{2}+\frac{\lambda}{2}\|\varphi\|^{2}-\left\langle\varphi, h-\mathrm{U}\left(b,t_{m}\right) \prod_{j=m}^{1}\left(
	\mathrm{I}+\mathrm{D}_{j}\right) \mathrm{U}\left(t_{j},t_{j-1}\right) x_{0}\right\rangle.
\end{align*}

The map $\varphi \rightarrow J_{\lambda}(\varphi)$ is continuous and strictly convex. The functional $J_{\lambda}(\cdot)$ admits a unique minimum $\widehat{\varphi}_{\lambda}$ that defines a map $\Phi: \mathbb{H} \rightarrow \mathbb{H}$. Since $J_{\lambda}(\varphi)$ is Frechet differentiable at $\widehat{\varphi}_{\lambda}$, by the optimality of $\widehat{\varphi}_{\lambda}$, we must have

\begin{align}\label{eqn2.4}
	\frac{d}{d \varphi} J_{\lambda}(\varphi)=&\Theta_{0}^{t_{m}} \widehat{\varphi}_{\lambda}+\Gamma_{t_{m}}^{b} \widehat{\varphi}_{\lambda}+\widetilde{\Theta}_{0}^{t_{m}} \widehat{\varphi}_{\lambda}+\widetilde{\Gamma}_{t_{m}}^{b} \widehat{\varphi}_{\lambda}+\lambda \widehat{\varphi}_{\lambda}-h \nonumber\\
	&\quad+\mathrm{U}\left(b,t_{m}\right) \prod_{j=m}^{1}\left(\mathrm{I}+\mathrm{D}_{j}\right) \mathrm{U}\left(t_{j},t_{j-1}\right) x_{0}=0.
\end{align}

By solving above equation \ref{eqn2.4}  for $\widehat{\varphi}_{\lambda}$, we get

\begin{align}\label{eqn2.5}
	\widehat{\varphi}_{\lambda}=\left(\lambda \mathrm{I}+\Theta_{0}^{t_{m}}+\Gamma_{t_{m}}^{b}+\widetilde{\Theta}_{0}^{t_{m}}+\widetilde{\Gamma}_{t_{m}}^{b}\right)^{-1}\left(h-\mathrm{U}\left(b,t_{m}\right) \prod_{j=m}^{1}\left(\mathrm{I}+\mathrm{D}_{j}\right) \mathrm{U}\left(t_{j},t_{j-1}\right) x_{0}\right). 
\end{align}
Now we define control $u(s)$ as following
\begin{align}\label{eqn2.6}
	\begin{split}
		u(s)=\bigg(\sum_{k=1}^{m}\mathrm{B}^{*}\mathrm{U}(t_{k},s)^{*}\prod_{i=k+1}^{m} \mathrm{U}(t_{i},t_{i-1})^{*}\mathrm{U}(b,t_{m})^{*}\chi(t_{k-1},t_{k})\\
		+\mathrm{B}^{*}\mathrm{U}(b,s)^{*}\chi(t_{m},b)\bigg)\widehat{\varphi}_{\lambda},\\
		v_{m}=\mathrm{E}_{m}^{*}\mathrm{U}(b,t_{m})^{*}\widehat{\varphi}_{\lambda},\quad
		v_{k}=\mathrm{E}_{k}^{*}\prod_{i=k}^{m} \mathrm{U}(t_{i},t_{i-1})^{*}(\mathrm{I}+\mathrm{D}_{i}^{*})\mathrm{U}(b,t_{m})^{*}\widehat{\varphi}_{\lambda}.
	\end{split}
\end{align}
Let $x_{\lambda}(b)$ be the solution at the final point $b$ corresponding to the above defined control, can be expressed as:
\begin{align*}
	x_{\lambda}(b)= \mathrm{U}(b,t_{m})\prod_{j=p}^{1} \left(\mathrm{I}+ \mathrm{D}_{j}\right)\left(\mathrm{U}(t_{j}-t_{j-1})\right)x_{0}+\Theta_{0}^{t_{m}} \widehat{\varphi}_{\lambda}+\Gamma_{t_{p}}^{b} \widehat{\varphi}_{\lambda}+\widetilde{\Theta}_{0}^{t_{p}} \widehat{\varphi}_{\lambda}+\widetilde{\Gamma}_{t_{p}}^{b} \widehat{\varphi}_{\lambda}+\lambda \widehat{\varphi}_{\lambda}.
\end{align*}.
Now from \ref{eqn2.3}, \ref{eqn2.5} and \ref{eqn2.6} we get
\begin{align}\label{2.7}
	x_{\lambda}(b)-h=-\lambda \left(\lambda I+\Theta_{0}^{t_{m}}+\Gamma_{t_{m}}^{b}+\widetilde{\Theta}_{0}^{t_{m}}+\widetilde{\Gamma}_{t_{m}}^{b}\right)^{-1}\left(h-\mathrm{U}\left(b,t_{m}\right) \prod_{j=m}^{1}\left(I+D_{j}\right) T\left(t_{j}-t_{j-1}\right) x_{0}\right), 
\end{align}
The above expression shows that the linear system \ref{4} is approximately controllable iff $-\lambda \left(\lambda I+\Theta_{0}^{t_{m}}+\Gamma_{t_{m}}^{b}+\widetilde{\Theta}_{0}^{t_{m}}+\widetilde{\Gamma}_{t_{m}}^{b}\right)^{-1}$ converges to zero operator as $\lambda \mapsto 0^{+}$ in strong operator topology. Therefore,  $(a)\Longleftrightarrow (d)$.

In order to establish the existence results for the system \eqref{1}, we require the following assumptions:
	\begin{Ass}
		\begin{enumerate}
			\item[(A1)] 	For every $x\in\mathbb{H}$, $\lambda\left(\lambda \mathrm{I}+\Theta_{0}^{t_{m}}+\Gamma_{t_{m}}^{b}+\widetilde{\Theta}_{0}^{t_{m}}+\widetilde{\Gamma}_{t_{m}}^{b}\right)^{-1}$ converges to zero operator as $\lambda \rightarrow 0^{+}$ in strong operator topology.
			
			\item[(A2)] \begin{enumerate} 
				\item [(i)]  The function $f:[0,b]\times \mathbb{H}\rightarrow \mathbb{H}$ is continuous and there is a constant $L_{f}$ such that for every $t \in[0, b]$ and $x, y \in \mathbb{H}$,
				\begin{align*}
					\|f(t, x)-f(t, y)\| \leq L_{f}\|x-y\|,
				\end{align*}
				
				\item [(ii)]  there exists $C_{f}$ such that $\left\|f\left(t,x(t)\right)\right\|\le C_{f}$ for $t\in [0,b]$.	
			\end{enumerate}
		\item[(A3)] \begin{enumerate} 
			\item [(i)]  The function $\xi :[0,b]\times \mathbb{H}\rightarrow \mathbb{H}$ is continuous and there is a constant $\tilde{L}_{\xi}$ such that for every $t \in[0, b]$ and $x, y \in \mathbb{H}$,
			\begin{align*}
				\|\xi (t, x)-\xi(t, y)\| \leq L_{\xi}\|x-y\|,
			\end{align*}
			
			\item [(ii)]  there exists $C_{\xi}$ such that $\left\|\xi\left(t,x(t)\right)\right\|\le C_{\xi}$ for $t\in [0,b]$.			
		\end{enumerate}
		\end{enumerate}
	
	\end{Ass}
	\begin{rem}
	Note that assumptions on $f$ and $\xi$ can be relaxed according to fixed point theorem which we are appying.
	\end{rem}
	
	\section{\textbf{Existence and Approximate Controllability of Semilinear System}}
	
The primary goal of this section is to identify sufficient conditions for the solvability of system \ref{1}. To achieve this, we will first demonstrate that, for each $\lambda$ and a fixed $h \in \mathbb{H}$, system \ref{1} possesses at least one mild solution.
We prove the existence of a mild solution of the system \eqref{1} with the control 
	\begin{align}\label{5}
		\begin{split}
			u(s)=\bigg(\sum_{k=1}^{m}\mathrm{B}^{*}\mathrm{U}(t_{k},s)^{*}\prod_{i=k+1}^{m} \mathrm{U}(t_{i},t_{i-1})^{*}\mathrm{U}(b,t_{m})^{*}\chi(t_{k-1},t_{k})\\
			+\mathrm{B}^{*}\mathrm{U}(b,s)^{*}\chi(t_{m},b)\bigg)\widehat{\varphi}_{\lambda},\\
			v_{m}=\mathrm{E}_{m}^{*}\mathrm{U}(b,t_{m})^{*}\widehat{\varphi}_{\lambda},\quad
			v_{k}=\mathrm{E}_{k}^{*}\prod_{i=k}^{m} \mathrm{U}(t_{i},t_{i-1})^{*}(\mathrm{I}+\mathrm{D}_{i}^{*})\mathrm{U}(b,t_{m})^{*}\widehat{\varphi}_{\lambda},
		\end{split}
	\end{align}
	with 
	\begin{align*}
		\widehat{\varphi}_{\lambda}=&\bigg(\lambda \mathcal{I}+\Theta_{0}^{t_{m}}+\Gamma_{t_{m}}^{b}+\widetilde{\Theta}_{0}^{t_{m}}+\widetilde{\Gamma}_{t_{m}}^{b}\bigg)^{-1} \times g\left(x(\cdot)\right),
	\end{align*} 
	where
	\begin{align*}
		g\left(x(\cdot)\right)=&\bigg(h-\mathrm{U}\left(b,t_{m}\right) \prod_{j=m}^{1}(\mathrm{I}+\mathrm{D}_{j}) \mathrm{U}\left(t_{j},t_{j-1}\right) x_{0}\\&\quad-\int_{t_{m}}^{b} \mathrm{U}(b,s)\left(f\left(s,x(s)\right)+\int_{0}^{s}q(s-\tau)\xi (\tau,x(\tau)) d\tau\right) d s\\&\qquad-\mathrm{U}\left(b,t_{m}\right)\sum_{i=1}^{m} \prod_{j=m}^{i+1}\left(\mathrm{I}+\mathrm{D}_{j}\right) \mathrm{U}\left(t_{j},t_{j-1}\right)\left(\mathrm{I}+\mathrm{D}_{i}\right)\\&\qquad \qquad\times\int_{t_{i-1}}^{t_{i}} \mathrm{U}\left(t_{i},s\right)\left(f\left(s,x(s)\right)+\int_{0}^{s}q(s-\tau)\xi (\tau,x(\tau)) d\tau\right) d s\bigg).
	\end{align*}
	
	With these assumptions established, we are now ready to prove the existence and uniqueness of the mild solution for \eqref{1} using the fixed point theorem \ref{fixed point thm}.

	\begin{thm} \label{Existence}
	    If the assumptions (R1)-(R4) and (A2)-(A3)
		are satisfied. Then for every $\lambda>0$ and for fixed $h\in \mathbb{H}$, the system \eqref{1}  has at least one mild solution in $\mathcal{PC}\left([0, b],\mathbb{H}\right)$ provided that
		\begin{equation}\label{C1}
			\max \left\{\mathcal{N}, \mathcal{K}_{1}\right\}<1, 
		\end{equation}
		and
		\begin{equation}\label{C2}
			\max \left\{M ; \mathscr{L}\right\}<1,
		\end{equation}
		where $\mathcal{N}$ and $\mathcal{K}_{1}$ are given by:
		$$
		\left\{\begin{array}{l}
			\mathcal{N}=M+\frac{ M^{3} M_{\mathrm{B}}^{2}b}{\lambda}, \\
			\mathcal{K}_{1}=M^{k+1}\prod_{j=i}^{k}\left(1+\left\|\mathrm{D}_{j}\right\|\right)\bigg(1+\frac{M^{2} M_{B}^{2} b}{\lambda}\left(mM^{m-k}+1\right)(M N+1)\\
			\hspace{5.5cm}+\mathcal{K}_{0}+\frac{M^{2}}{\lambda}\|\mathrm{E}_{k}\|\left\|\mathrm{E}_{k}^{*}\prod_{i=k}^{m}\left\|\mathrm{T}^{*}(t_i-t_{i-1})\right\|\left(\mathcal{I}+\mathrm{D}_{i}^{*}\right)\right\|\bigg),\\
			\mathcal{K}_{0}=\frac{M^{2}}{\lambda}\sum_{i=2}^{k} \prod_{j=i}^{k}\left(1+\left\|\mathrm{D}_{j}\right\|\right)\left\|\mathrm{T}(t_j-t_{j-1})\right\|\left\|\mathrm{E}_{i-1}\right\|\left\|\mathrm{E}_{i-1}^{*}\prod_{l=i-1}^{m}\left\|\mathrm{U}(t_l,t_{l-1})^{*}\right\|\left(\mathcal{I}+\mathrm{D}_{l}^{*}\right)\right\|,\\
			C_{i}=\prod_{j=k}^{i+1}\left(1+\left\|D_{j}\right\|\right)\left\|T\left(t_{j}-t_{j-1}\right)\right\|\left(1+\left\|D_{i}\right\|\right), \quad N=\sum_{i=1}^{k} C_{i},\\
			\mathscr{L}= M^{k+1} \prod_{j=1}^{k}\left(1+\left\|\mathrm{D}_{j}\right\|\right)+M^{2} Nb \left(L_{f}+q^{*} L_{\xi}\right) \text{ and } q^{*}=\int_{0}^{t}|q(t-s)| ds.
			
		\end{array}\right.
		$$
	\end{thm} 
	Proof. For each constant $r_{0}>0$, let
	
	\begin{align*}
		\mathcal{B}_{r_{0}}=\left\{x \in \mathcal{PC}\left([0, b],\mathbb{H}\right): \quad\|x\|_{\mathcal{PC}} \leq r_{0}\right\}.
	\end{align*}
	
	It is easy to see that $\mathcal{B}_{r_{0}}$ is a bounded closed convex set. Define operators $F_{1}$ and $F_{2}$ on $\mathcal{B}_{r_{0}}$ as follows:
	
	$$
	\begin{gathered}
		\left(F_{1} x\right)(t)=\left\{\begin{array}{l}
			\mathrm{U}(t,0) x_{0}, \quad \text { for } \quad t_{0}<t \leq t_{1} \\
			\mathrm{U}\left(t,t_{k}\right) \prod_{j=k}^{1}\left(\mathcal{I}+\mathrm{D}_{j}\right) \mathrm{U}\left(t_{j},t_{j-1}\right) x_{0} \\
			+\mathrm{U}\left(t,t_{k}\right) \sum_{i=1}^{k} \prod_{j=k}^{i+1}\left(\mathcal{I}+\mathrm{D}_{j}\right) \mathrm{U}\left(t_{j},t_{j-1}\right)\left(\mathcal{I}+\mathrm{D}_{i}\right) \int_{t_{i-1}}^{t_{i}} \mathrm{U}\left(t_{i},s\right) \mathrm{B} u(s) d s \\
			+\mathrm{U}\left(t,t_{k}\right) \sum_{i=1}^{k} \prod_{j=k}^{i+1}\left(\mathcal{I}+\mathrm{D}_{j}\right) \mathrm{U}\left(t_{j},t_{j-1}\right)\left(\mathcal{I}+\mathrm{D}_{i}\right)\\ \qquad\times\int_{t_{i-1}}^{t_{i}} \mathrm{U}\left(t_{i},s\right) \left[f\left(s,x(s)\right)+\int_{0}^{s}q(s-\tau)\xi (\tau,x(\tau)) d\tau\right] d s\\
			+\mathrm{U}\left(t,t_{k}\right) \sum_{i=2}^{k} \prod_{j=k}^{i}\left(\mathcal{I}+\mathrm{D}_{j}\right) \mathrm{U}\left(t_{j},t_{j-1}\right) \mathrm{E}_{i-1} v_{i-1}+\mathrm{U}\left(t,t_{k}\right) \mathrm{E}_{k} v_{k}, \\
			\text { for } t_{k}<t \leq t_{k+1}, k \geq 1,
		\end{array}\right. \\
		\left(F_{2} x\right)(t)=\left\{\begin{array}{l}
			\int_{0}^{t} \mathrm{U}(t,s)\left[\mathrm{B} u(s)+f\left(s,x(s)\right)+\int_{0}^{s}q(s-\tau)\xi (\tau,x(\tau)) d\tau\right] d s,  \text { for }  t_{0}<t \leq t_{1} \\
			\int_{t_{k}}^{t} \mathrm{U}(t,s)\left[\mathrm{B} u(s)+f\left(s,x(s)\right)+\int_{0}^{s}q(s-\tau)\xi (\tau,x(\tau)) d\tau\right] d s, \text { for } t_{k}<t \leq t_{k+1}, k \geq 1.
		\end{array}\right.
	\end{gathered}
	$$
	Clearly, $x$ is a mild solution of (2) if and only if the operator equation $x=F_{1} x+F_{2} x$ has a solution. To establish this, we will demonstrate that the operator $F_{1}+F_{2}$ has\\
	a fixed point by applying Theorem \ref{fixed point thm} . For this, we proceed in several steps.
	
	\textbf{Step 1:}To prove that there exists a positive number $r_{0}$ such that $F_{1} x+F_{2} y \in \mathcal{B}_{r_{0}}$ whenever $x, y \in \mathcal{B}_{r_{0}}$, we proceed as follows:
	
	Choose
	
	\begin{equation*}
		r_{0} \geq \max \left(\frac{\left(\frac{ M^{2} M_{\mathrm{B}}^{2}b}{\lambda}\|h\|+\frac{ M^{2} M_{\mathrm{B}}^{2} b}{\lambda}\left(M b C_{f}+ M b C_{\xi }q^{*}\right)+M b C_{f}+ M b C_{\xi }q^{*}\right)}{1-\mathcal{N}},\frac{\mathcal{K}_{2}}{1-\mathcal{K}_{1}}\right),
	\end{equation*}
	where 
	\begin{align*}
		\mathcal{K}_2&=\left(\frac{M^{2}M_{\mathrm{B}}^{2}b}{\lambda}\left(MN+1\right)\left(mM^{m-k}+1\right)+\mathcal{K}_{0}\frac{M^{2}}{\lambda}\|\mathrm{E}_{k}\|\left\|\mathrm{E}_{k}^{*}\prod_{i=k}^{m}\left\|\mathrm{U}(t_i,t_{i-1})^{*}\right\|\left(\mathcal{I}+\mathrm{D}_{i}^{*}\right)\right\|\right) \\& \qquad \times\left(\|h\|+M \left(C_{f}+q^{*} C_{\xi}\right) b + M^{2}N C b\right)+ M b \left(C_{f}+q^{*} C_{\xi}\right)(MN+1).
	\end{align*}
	
	First, we calculate for $t_{0}<t \leq t_{1}$ and $s\in [0,b]$,
	\begin{align*}
		u(s)=\mathrm{B}^{*}\mathrm{U}x(t_1,s)^{*}\big(\lambda I+\Gamma_{0}^{t_1}\big)^{-1}\left[h-\mathrm{U}(t,0)x_0-\int_{0}^{t_1}\mathrm{U}(t_1,s)\left(f\left(s,x(s)\right)+\int_{0}^{s}q(s-\tau)\xi (\tau,x(\tau)) d\tau\right)d s\right].
	\end{align*}
	Using the triangle inequality, Lipschitz conditions, and the boundedness of the evolution family $\mathrm{U}(t,s)$, the  norm $	\|u(s)\|_{\mathbb{U}}$ can be calculated as:
	\begin{align*}
		\|u(s)\|_{\mathbb{U}}&=\bigg\|\mathrm{B}^{*}\mathrm{U}(t_1,s)^{*}\big(\lambda I+\Gamma_{0}^{t_1}\big)^{-1}\\&\qquad\times\left[h-\mathrm{U}(t,0)x_0-\int_{0}^{t_1}\mathrm{U}(t_1,s)\left(f\left(s,x(s)\right)+\int_{0}^{s}q(s-\tau)\xi (\tau,x(\tau)) d\tau\right) d s\right]\bigg\|,\\
		&\le \left\|\mathrm{B}^{*}\right\|_\mathcal{L}\left\|\mathrm{U}(t_1-s)^{*}\right\|_{\mathbb{H}}\left\|\big(\lambda I+\Gamma_{0}^{t_1}\big)^{-1}\right\|\\&\qquad\times\bigg\|h-\mathrm{U}(t,0)x_0-\int_{0}^{t_1}\mathrm{U}(t_1,s)\left(f\left(s,x(s)\right)+\int_{0}^{s}q(s-\tau)\xi (\tau,x(\tau)) d\tau\right) d s\bigg\|,\\
		&\le \frac{ M M_{B}}{\lambda}\bigg(\|h\|+\left\|\mathrm{U}(t_1,0)\right\|_{\mathbb{H}}\|x_0\|\\&\qquad+\left\|\mathrm{U}(t_1,s)\right\|_{\mathbb{H}}\int_{0}^{t_1}\left\|\left(f\left(s,x(s)\right)+\int_{0}^{s}q(s-\tau)\xi (\tau,x(\tau)) d\tau\right)\right\|d s\bigg),\\
		&\le \frac{ M M_{B}}{\lambda}\left(\|h\|+Mr_0+M C_{f} b+ M C_{\xi}q^{*}\right).
	\end{align*}
	To calculate the norm of $u$ for $t_k<t\le t_{k+1}$, $k\ge1$ and $s\in [0,b]$ , first we find the norm of $\widehat{\varphi}_{\lambda}$ as follows:
	\begin{align*}
		\left \|\widehat{\varphi}_{\lambda}\right \|&\le \left\|\bigg(\lambda \mathcal{I}+\Theta_{0}^{t_{m}}+\Gamma_{t_{m}}^{b}+\widetilde{\Theta}_{0}^{t_{m}}+\widetilde{\Gamma}_{t_{m}}^{b}\bigg)^{-1}\right\|\\
		&\quad\times\Bigg(\left\|h\right\|+\big\|\mathrm{U}\left(b,t_{m}\right) \prod_{j=m}^{1}(\mathrm{I}+\mathrm{D}_{j}) \mathrm{U}\left(t_{j},t_{j-1}\right) x_{0}\big\|\\& \qquad \qquad+\left\|\int_{t_{m}}^{b} \mathrm{U}(b,s)\left(f\left(s,x(s)\right)+\int_{0}^{s}q(s-\tau)\xi (\tau,x(\tau)) d\tau\right) d s\right\|\\
		&\qquad \qquad+\bigg\|\mathrm{U}\left(b,t_{m}\right)\sum_{i=1}^{m} \prod_{j=m}^{i+1}\left(\mathrm{I}+\mathrm{D}_{j}\right) \mathrm{U}\left(t_{j}-t_{j-1}\right)\left(\mathrm{I}+\mathrm{D}_{i}\right)\\& \qquad \qquad \qquad\times\int_{t_{i-1}}^{t_{i}} \mathrm{U}\left(t_{i}-s\right)\left(f\left(s,x(s)\right)+\int_{0}^{s}q(s-\tau)\xi (\tau,x(\tau)) d\tau\right) d s\bigg\|\Bigg),\\
		&\le \frac{1}{\lambda}\Bigg(\left\|h\right\|+M^{k+1}\prod_{j=1}^{k}\left(1+\left\|\mathrm{D}_{j}\right\|\right)\left\|x_{0}\right\|+M b (C_{f}+ q^{*} C_{\xi})\\&\qquad \quad+M^{2}\sum_{i=1}^{m} C_{i}\int_{t_{i-1}}^{t_{i}} \left\|\left(f\left(s,x(s)\right)+\int_{0}^{s}q(s-\tau)\xi (\tau,x(\tau)) d\tau\right)\right\| d s\Bigg),\\
		&\le \frac{1}{\lambda}\Bigg(\left\|h\right\|+M^{k+1}\prod_{j=1}^{k}\left(1+\left\|\mathrm{D}_{j}\right\|\right)\left\|x_{0}\right\|+M b (C_{f}+ q^{*} C_{\xi})+M^{2} N b (C_{f}+ q^{*} C_{\xi})\Bigg).
	\end{align*}
	With the above help we can find the norm of $u$ as follows:
	\begin{align*}
		\left\|u(s)\right\|_{\mathbb{U}}&\le\\ &\left\|\bigg(\sum_{k=1}^{m}\mathrm{B}^{*}\mathrm{U}(t_{k},s)^{*}\prod_{i=k+1}^{m} \mathrm{U}^{*}(t_{i},t_{i-1})\mathrm{U}^{*}(b,t_{m})\chi(t_{k-1},t_{k})
		+\mathrm{B}^{*}\mathrm{U}(b,s)^{*}\chi(t_{m},b)\bigg)\right\|
		\left\|\widehat{\varphi}_{\lambda}\right\|,\\
		&\le \frac{1}{\lambda}\left(mM_{\mathrm{B}}M^{m+1-k}+M_{\mathrm{B}}M\right)\times\\
		&\qquad \Bigg(\left\|h\right\|+M^{k+1}\prod_{j=1}^{k}\left(1+\left\|\mathrm{D}_{j}\right\|\right)r_0+M b (C_{f}+ q^{*} C_{\xi}) + M^{2} N(C_{f}+ q^{*} C_{\xi})\Bigg).
	\end{align*}
	Now, for $0\le t\le  t_1$
	\begin{align*}
		\left\|\left(F_{1} x\right)(t)+\left(F_{2} x\right)(t)\right\| &\leq \|\mathrm{U}(t,0) x(0)\|_{\mathbb{H}}+\left\|\int_{0}^{t} \mathrm{U}(t,s) \mathrm{B} u(s) d s\right\|\\& \qquad +\left\|\int_{0}^{t} \mathrm{U}(t,s) \left(f\left(s,x(s)\right)+\int_{0}^{s}q(s-\tau)\xi (\tau,x(\tau)) d\tau\right) d s\right\|,\\
		& \leq  M\left\|x_{0}\right\|+ M\|\mathrm{B}\|_\mathcal{L} \int_{0}^{t}u(s) d s+ M b C_{f}+ M b C_{\xi} q^{*} , \\
		& \leq M r_{0}+\frac{ M^{2} M_{\mathrm{B}}^{2}b}{\lambda}\left(\|h\|+Mr_0+M b C_{f}+ M b C_{\xi} q^{*}\right) + M b C_{f}+ M b C_{\xi} q^{*}, \\
		& =\mathcal{N} r_{0}+\left(\frac{ M^{2} M_{\mathrm{B}}^{2}b}{\lambda}\|h\|+\left(\frac{ M^{2} M_{\mathrm{B}}^{2} b}{\lambda}+1\right)\left(M b C_{f}+ M b C_{\xi }q^{*}\right)\right),\\
		& \leq r_{0}.
	\end{align*}
	
	%
	
	For $t_{k}<t \leq t_{k+1}$ for $k \geq 1$, we have,

	\begin{align*}
		\left\|\left(F_{1} x\right)(t)+\left(F_{2} x\right)(t)\right\| &\leq  \left\|\mathrm{U}\left(t,t_{k}\right) \prod_{j=k}^{1}\left(\mathcal{I}+\mathrm{D}_{j}\right) T\left(t_{j}-t_{j-1}\right) x_{0}\right\| \\
		& \quad+\left\|\mathrm{U}\left(t,t_{k}\right) \sum_{i=1}^{k} \prod_{j=k}^{i+1}\left(\mathcal{I}+\mathrm{D}_{j}\right) \mathrm{U}\left(t_{j},t_{j-1}\right)\left(\mathcal{I}+\mathrm{D}_{i}\right) \int_{t_{i-1}}^{t_{i}} \mathrm{U}\left(t_{i},s\right) \mathrm{B} u(s) d s\right\|\\
		&\quad +\bigg\|\mathrm{U}\left(t,t_{k}\right) \sum_{i=1}^{k} \prod_{j=k}^{i+1}\left(\mathcal{I}+\mathrm{D}_{j}\right) \mathrm{U}\left(t_{j},t_{j-1}\right)\left(\mathcal{I}+\mathrm{D}_{i}\right)\\&\qquad \quad \times \int_{t_{i-1}}^{t_{i}} \mathrm{U}\left(t_{i},s\right) \left(f\left(s,x(s)\right)+\int_{0}^{s}q(s-\tau)\xi (\tau,x(\tau)) d\tau\right) d s\bigg\| \\
		&\quad +\left\|\mathrm{U}\left(t,t_{k}\right) \sum_{i=2}^{k} \prod_{j=k}^{i}\left(\mathcal{I}+\mathrm{D}_{j}\right) \mathrm{U}\left(t_{j},t_{j-1}\right) \mathrm{E}_{i-1} v_{i-1}\right\| \\
		&\quad + \left\|\mathrm{U}\left(t,t_{k}\right) \mathrm{E}_{k} v_{k}\right\|+\left\|\int_{t_{k}}^{t} \mathrm{U}(t,s) \mathrm{B} u(s) d s\right\|\\ & \quad +\left\|\int_{t_{k}}^{t} \mathrm{U}(t,s) \left(f\left(s,x(s)\right)+\int_{0}^{s}q(s-\tau)\xi (\tau,x(\tau)) d\tau\right) d s\right\|,\\
		& \leq  M^{ k+1} \prod_{j=1}^{k}\left(1+\left\|\mathrm{D}_{j}\right\|\right) r_{0}+M^{2}\|\mathrm{B}\|_\mathcal{L} \sum_{i=1}^{k} C_{i}  \int_{t_{i-1}}^{t_{i}} \|u(s)\|\hspace{.2cm} d s\\
		& \quad + M^{2} \sum_{i=1}^{k} C_{i}\int_{t_{i-1}}^{t_{i}}\|\left(f\left(s,x(s)\right)+\int_{0}^{s}q(s-\tau)\xi (\tau,x(\tau)) d\tau\right)\| d s\\& \quad+ M \sum_{i=2}^{k} \prod_{j=i}^{k}\left(1+\left\|\mathrm{D}_{j}\right\|\right)\left\|\mathrm{U}(t_j,t_{j-1})\right\|\left\|\mathrm{E}_{i-1}\right\| \left\|v_{i-1}\right\| +M\left\|\mathrm{E}_{k}\right\|\left\|v_{k}\right\|\\& \quad+M\|\mathrm{B}\|_\mathcal{L} \int_{t_k}^{t}u(s) d s+M\int_{t_{k}}^{t}\|\left(f\left(s,x(s)\right)+\int_{0}^{s}q(s-\tau)\xi (\tau,x(\tau)) d\tau\right)\| d s,\\
		&\leq \mathcal{K}_{1}r_0+ \mathcal{K}_{2},\\
		&\leq r_{0}.	
	\end{align*}
	

	Consequently, $F_{1}+F_{2}$ maps $\mathcal{B}_{r_{0}}$ to $\mathcal{B}_{r_{0}}$.
	
	\textbf{Step 2:} The next step is to prove that $F_{1}$ is a contraction.
	
	To demonstrate that $F_{1}$ is a contraction mapping on the set $\mathcal{B}_{r}$, it is necessary to show that there exists a constant $0 < \mathscr{L} < 1$ such that for all $x, y \in \mathcal{B}_{r}$,
		\begin{align*}
		\left\|F_{1} x-F_{1} y\right\|_{\mathcal{PC}} \leq \mathscr{L}\|x-y\|_{\mathcal{PC}}.
	\end{align*}
	
	Let $x, y \in \mathcal{B}_{r}$. We will estimate $\left\|F_{1} x-F_{1} y\right\|_{\mathcal{PC}}$ for $t_{0}<t \leq t_{1}$ and $t_{k}<t \leq t_{k+1}$.
	
	For $t_{0}<t \leq t_{1}$ :
	
	\begin{align*}
		\left\|\left(F_{1} x\right)(t)-\left(F_{1} y\right)(t)\right\|=\|\mathrm{U}(t,0)(x(0)-y(0))\|_{\mathbb{H}}.
	\end{align*}
	
	Using the properties of the evolution operator $\mathrm{U}(t,s)$ :
	
	\begin{align*}
		\|\mathrm{U}(t,0)(x(0)-y(0))\| \leq M \|x(0)-y(0)\|,
	\end{align*}
	for $M<1$, $F_{1}$ is a contraction map.

	For $t_{k}<t \leq t_{k+1}, k \geq 1$ :
	\begin{align*}
		& \left\|\left(F_{1} x\right)(t)-\left(F_{1} y\right)(t)\right\|^ \leq  \left\|\mathrm{U}\left(t,t_{k}\right) \prod_{j=k}^{1}\left(\mathcal{I}+\mathrm{D}_{j}\right) \mathrm{U}\left(t_{j},t_{j-1}\right)\left(x_{0}-y_{0}\right)\right\| \\
		& +\bigg\|\mathrm{U}\left(t,t_{k}\right) \sum_{i=1}^{k} \prod_{j=k}^{i+1}\left(\mathcal{I}+\mathrm{D}_{j}\right) \mathrm{U}\left(t_{j},t_{j-1}\right)\left(\mathcal{I}+\mathrm{D}_{i}\right) \int_{t_{i-1}}^{t_{i}} \mathrm{U}\left(t_{i},s\right)\\&\qquad\times\left[\left(f\left(s,x(s)\right)-f\left(s,y(s)\right)\right)+\int_{0}^{s}q(s-\tau)\left(\xi (\tau,x(\tau))-\xi (\tau,y(\tau))\right) d\tau \right] d s\bigg\|. 
	\end{align*}

	Using the properties of $\mathrm{U}(t,s)$, the boundedness of operators $\mathrm{D}_{j}$, and assumptions on $f$ :
	
	\begin{align*}
		\left\|\mathrm{U}\left(t,t_{k}\right) \prod_{j=k}^{1}\left(\mathcal{I}+\mathrm{D}_{j}\right) \mathrm{U}\left(t_{j},t_{j-1}\right)\left(x_{0}-y_{0}\right)\right\| \leq M^{k+1} \prod_{j=1}^{k}\left(1+\left\|\mathrm{D}_{j}\right\|\right)\left\|x_{0}-y_{0}\right\|.
	\end{align*}
	
	Since $x, y \in \mathcal{B}_{r}$ :
	\begin{align*}
		\left\|x_{0}-y_{0}\right\| \leq\|x-y\|_{\mathcal{PC}}.
	\end{align*}
	Thus,
	\begin{align*}
		\left\|\mathrm{U}\left(t,t_{k}\right) \prod_{j=k}^{1}\left(\mathcal{I}+\mathrm{D}_{j}\right) \mathrm{U}\left(t_{j},t_{j-1}\right)\left(x_{0}-y_{0}\right)\right\| \leq M^{k+1} \prod_{j=1}^{k}\left(1+\left\|\mathrm{D}_{j}\right\|\right)\|x-y\|_{\mathcal{PC}}.
	\end{align*}
	
	For the second term, using the properties of $\mathrm{U}(t,s)$ and $\mathrm{D}_{j}$, and assumption (A2) of $f$ :
	
	\begin{align*}
		& \bigg\|\mathrm{U}\left(t,t_{k}\right) \sum_{i=1}^{k} \prod_{j=k}^{i+1}\left(\mathcal{I}+\mathrm{D}_{j}\right) \mathrm{U}\left(t_{j},t_{j-1}\right)\left(\mathcal{I}+\mathrm{D}_{i}\right) \int_{t_{i-1}}^{t_{i}} \mathrm{U}\left(t_{i},s\right)\\&\quad \times \left[\left(f\left(s,x(s)\right)-f\left(s,y(s)\right)\right)+\int_{0}^{s}q(s-\tau)\left(\xi (\tau,x(\tau))-\xi (\tau,y(\tau))\right) d\tau \right] d s\bigg\|, \\
		& \leq M^{2} \Bigg(\sum_{i=1}^{k} \prod_{j=i+1}^{k}\left(1+\left\|\mathrm{D}_{j}\right\|\right)\left\|\mathrm{U}\left(t_{j},t_{j-1}\right)\right\|_{\mathbb{H}}\left(1+\left\|\mathrm{D}_{i}\right\|\right)\\&\quad\times \int_{t_{i-1}}^{t_{i}}\big\|\left(f\left(s,x(s)\right)-f\left(s,y(s)\right)\right)+\int_{0}^{s}q(s-\tau)\left(\xi (\tau,x(\tau))-\xi (\tau,y(\tau))\right) d\tau\big\| d s\Bigg), \\
		& \leq M^{2}\bigg(\sum_{i=1}^{k} C_{i} \int_{t_{i-1}}^{t_{i}}\left[\left\|(f\left(s,x(s)\right)-f\left(s,y(s)\right))\right\|+\int_{0}^{s}\|q(s-\tau)\|\left\|\xi (\tau,x(\tau))-\xi (\tau,y(\tau))\right\| d\tau\right] d s\bigg), \\
		& \leq M^{2} N b \left(L_{f}+q^{*}L_{\xi}\right) \left\|x-y\right\|_{\mathcal{PC}}.
	\end{align*}
	
	Combining all terms, we get:\\
	$\left\|\left(F_{1} x\right)(t)-\left(F_{1} y\right)(t)\right\| \leq\left( M^{k+1} \prod_{j=1}^{k}\left(1+\left\|\mathrm{D}_{j}\right\|\right)+M^{2} Nb \left(L_{f}+q^{*}L_{\xi}\right)\right)\|x-y\|_{\mathcal{PC}}$.
	
	To show that $F_{1}$ is a contraction, we need the right-hand side to be less than $\|x-y\|_{\mathcal{PC}}$. Hence, we need
	\begin{align*}
		M^{k+1} \prod_{j=1}^{k}\left(1+\left\|\mathrm{D}_{j}\right\|\right)+M^{2} Nb \left(L_{f}+q^{*}L_{\xi}\right)<1.
	\end{align*}
	
Therefore, there exists a constant $\mathscr{L}\in (0,1)$ such that:
	\begin{align*}
		\left\|F_{1} x-F_{1} y\right\|_{\mathcal{PC}} \leq \mathscr{L}\|x-y\|_{\mathcal{PC}}.
	\end{align*}
	
	This shows that on $\mathcal{B}_{r_{0}}$, $F_{1}$ is a contraction map .
	
	\textbf{Step 3:} Now we will show that $F_{2} $ is continuous and $F_{2}(\mathcal{B}_{r_{0}})$ is relatively compact subset of $\mathcal{B}_{r_{0}}$.
	
	First, we need to prove that the mapping $F_{2}$ is continuous on $\mathcal{B}_{r_{0}}$. To do this, let $x_{n} \rightarrow x$ in $\mathcal{B}_{r_{0}}$. Then, we have:
	
	\begin{align*}
		f\left(t, x_{n}(t)\right) \rightarrow f\left(t, x(t)\right), \text{ and } \xi\left(t, x_{n}(t)\right) \rightarrow \xi\left(t, x(t)\right)\text { as } n \rightarrow \infty.
	\end{align*}
	Moreover, for $t_{0} \leq t \leq t_{1}$, by Lebesgue dominated convergence theorem, we can get
	
	\begin{align*}
		& \left\|\int_{0}^{t} \mathrm{U}(t,s)\left[f\left(t, x_{n}(t)\right)+\int_{0}^{s} q(s-\tau)\xi(\tau,x_{n}(\tau)) d\tau- f\left(t, x(t)\right)-\int_{0}^{s} q(s-\tau)\xi(\tau,x(\tau)) d\tau\right] d s\right\| \\
		& \quad\leq M \int_{0}^{t}\left[\left\|f\left(t, x_{n}(t)\right)-f\left(t, x(t)\right)\right\|+\int_{0}^{s}\| q(s-\tau)\| \left\|\xi(\tau,x_{n}(\tau))-\xi(\tau,x(\tau))\right\|\right] d s \rightarrow 0,\\& \qquad \text { as } n \rightarrow \infty.
	\end{align*}
	\begin{align*}
		\left\|F_{2}\left(x_{n}\right)-F_{2}(x)\right\| &\leq \bigg\|\int_{0}^{t} \mathrm{U}(t,s)\big[f\left(t, x_{n}(t)\right)+\int_{0}^{s} q(s-\tau)\xi(\tau,x_{n}(\tau)) d\tau\\& \qquad\qquad \qquad \qquad- f\left(t, x(t)\right)-\int_{0}^{s} q(s-\tau)\xi(\tau,x(\tau)) d\tau\big] d s\bigg\| \\
		& \rightarrow 0 \text { as } n \rightarrow \infty.
	\end{align*}
For $t_{k} < t \leq t_{k+1}$ with $k \geq 1$, the argument is similar to that for $t_{0} < t \leq t_{1}$. Hence, it follows that $F_{2}$ is continuous on $\mathcal{B}_{r_{0}}$.

Next, we demonstrate that for any $t \in [0, b]$, the set $\mathscr{V}(t) = \{F_{2}(x)(t) \mid x \in \mathcal{B}_{r_{0}}\}$ is relatively compact in $\mathbb{H}$. To establish this, we will utilize the extended version of the Ascoli-Arzelà theorem (Theorem 2.1, \cite{wei2006nonlinear}). For $t = 0$, it is evident that $\mathscr{V}(0)$ is relatively compact in $\mathbb{H}$. Now, for $0 < t \leq b$, let $\epsilon \in (0, t)$. By applying Lemma \ref{lem2.2}, we find that the operator $\mathrm{U}(t, t - \epsilon)$ is compact. We define an operator $F^{\epsilon}$ on $\mathcal{B}_{r_{0}}$ by:
	$$
	\left(F^{\epsilon} x\right)(t)=\left\{\begin{array}{l}
		\int_{0}^{t-\epsilon} \mathrm{U}(t,s)\left[\mathrm{B} u(s)+f\left(s,x(s)\right)+\int_{0}^{s}q(s-\tau)\xi (\tau,x(\tau)) d\tau\right] d s \\
		=\mathrm{U}(t,t-\epsilon) \int_{0}^{t-\epsilon} \mathrm{U}(t-\epsilon,s)\\ \qquad \quad \times\left[\mathrm{B} u(s)+f\left(s,x(s)\right)+\int_{0}^{s}q(s-\tau)\xi (\tau,x(\tau)) d\tau\right] d s \text { if } t_{0}<t \leq t_{1}, \\
		\int_{t_{k}}^{t-\epsilon} \mathrm{U}(t,s)\left[\mathrm{B} u(s)+f\left(s,x(s)\right)+\int_{0}^{s}q(s-\tau)\xi (\tau,x(\tau)) d\tau\right] d s \\
		=\mathrm{U}(t,t-\epsilon) \int_{t_{k}}^{t-\epsilon} \mathrm{U}(t-\epsilon,s)\\ \qquad \quad \times\left[\mathrm{B} u(s)+f\left(s,x(s)\right)+\int_{0}^{s}q(s-\tau)\xi (\tau,x(\tau)) d\tau\right] d s \text { if } t_{k}<t \leq t_{k+1}, k \geq 1.
	\end{array}\right.
	$$
	
	Then the set $\left\{\left(F^{\epsilon}\right)(t): x \in \mathcal{B}_{r_0}\right\}$ is relatively compact in $\mathbb{H}$ because $\mathrm{U}(t,t-\epsilon)$ is compact. This compactness helps us establish the desired continuity properties. Now, let's consider the case for $t_{0}<t \leq t_{1}$ :
	
	\begin{align*}
		\left\|\left(F_{2} x\right)(t)-\left(F^{\epsilon} x\right)(t)\right\| & \leq  \left\|\int_{t-\epsilon}^{t} \mathrm{U}(t,s) \mathrm{B} u(s) d s\right\| \\
		&\quad+ \left\|\int_{t-\epsilon}^{t} \mathrm{U}(t,s) \left[f\left(s,x(s)\right)+\int_{0}^{s}q(s-\tau)\xi (\tau,x(\tau)) d\tau\right] d s\right\|.
	\end{align*}
	
	To estimate the component involving $\mathrm{B} u^{\lambda}(s)$, we apply the triangle inequality followed by the Cauchy-Schwarz inequality. This yields:
	
	\begin{align*}
		\left\|\int_{t-\epsilon}^{t} \mathrm{U}(t,s) \mathrm{B} u(s) d s\right\| \leq M M_{B} {\epsilon}^{\frac{1}{2}} \left(\int_{t-\epsilon}^{t}\|u(s)\|^{2} ds\right)^{\frac{1}{2}}
	\end{align*}
	Using assumptions (A2) and (A3), we have
	
	\begin{align*}
		\bigg\|\int_{t-\epsilon}^{t}& \mathrm{U}(t,s)\left[f\left(s,x(s)\right)+\int_{0}^{s}q(s-\tau)\xi (\tau,x(\tau)) d\tau\right] d s\bigg\|\\& \leq \left(\int_{t-\epsilon}^{t}\left\|\mathrm{U}(t,s) \left[f\left(s,x(s)\right)+\int_{0}^{s}q(s-\tau)\xi (\tau,x(\tau)) d\tau\right]\right\| d s\right), \\
		& \leq M \int_{t-\epsilon}^{t} \left\|\left(f\left(s,x(s)\right)+\int_{0}^{s}q(s-\tau)\xi (\tau,x(\tau)) d\tau\right)\right\| d s\\
		& \leq M \left(C_{f}-q^{*}C_{\xi}\right) \epsilon.
	\end{align*}

	Combining all terms, we get:
	
	$\left\|\left(F_{2} x\right)(t)-\left(F^{\epsilon} x\right)(t)\right\| \leq  M \left(C_{f}-q^{*}C_{\xi}\right) \epsilon+ M M_{B} {\epsilon}^{\frac{1}{2}} \left(\int_{t-\epsilon}^{t}\|u(s)\|^{2} ds\right)^{\frac{1}{2}}$.
	
	As $\epsilon \rightarrow 0$ :
	\begin{align*}
		\left\|\left(F_{2} x\right)(t)-\left(F^{\epsilon} x\right)(t)\right\| \rightarrow 0.
	\end{align*}
	
For $t_{k} < t \leq t_{k+1}$, with $k \geq 1$, the definitions of $F_{2}$ and $F^{\epsilon}$ allow us to derive similar results as previously discussed.

Therefore, since $F_{2} x$ can be approximated arbitrarily closely by $F^{\epsilon} x$, and $F^{\epsilon} x$ is relatively compact in $\mathbb{H}$, it follows that $\mathscr{V}(t) = \{F_{2}(x)(t) \mid x \in \mathcal{B}_{r_{0}}\}$ is relatively compact in $\mathbb{H}$.

Finally, we show that $F_{2}(\mathcal{B}_{r_{0}})$ is equicontinuous on $[0, b]$.

	Let $0\leq s_1\leq s_2 \leq t_1$ for any $x\in \mathcal{B}_{r_0}$, we consider the following estimate 
	\begin{align}\label{eq3.6}
		\big\|	F_{2}x(s_2)&-F_{2}x(s_1)\big\|\nonumber\\&\leq\left\|\int_{0}^{s_1}\left[\mathrm{U}(s_2,s)-\mathrm{U}(s_1,s)\right]\left[B u(s)+f\left(s,x(s)\right)+\int_{0}^{s}q(s-\tau)\xi (\tau,x(\tau)) d\tau\right]ds\right\|\nonumber\\
		&\quad + \left\|\int_{s_1}^{s_2}\mathrm{U}(s_2,s)\left[B u(s)+f\left(s,x(s)\right)+\int_{0}^{s}q(s-\tau)\xi (\tau,x(\tau)) d\tau\right]ds\right\|,\nonumber\\
		& \leq \int_{s_1}^{s_2}\left\|\mathrm{U}(s_2,s)\right\|_{\mathcal{L}(\mathbb{H})}\|\mathrm{B}\|_{\mathcal{L}} \|u(s)\|_{\mathbb{U}} ds\nonumber\\
		&\quad +\int_{s_1}^{s_2}\left\|\mathrm{U}(s_2,s)\right\|_{\mathcal{L}(\mathbb{H})} \left\|\left(f\left(s,x(s)\right)+\int_{0}^{s}q(s-\tau)\xi (\tau,x(\tau)) d\tau\right)\right\|\nonumber\\
		& \quad +\int_{0}^{s_1}\left\|\mathrm{U}(s_2,s)-\mathrm{U}(s_1,s)\right\|_{\mathcal{L}(\mathbb{H})}\|\mathrm{B}\|_{\mathcal{L}} \|u(s)\|_{\mathbb{U}} ds\nonumber\\
		&\quad +\int_{0}^{s_1}\left\|\mathrm{U}(s_2,s)-\mathrm{U}(s_1,s)\right\|_{\mathcal{L}(\mathbb{H})} \left\|\left(f\left(s,x(s)\right)+\int_{0}^{s}q(s-\tau)\xi (\tau,x(\tau)) d\tau\right)\right\| ds \nonumber\\
		&\leq M M_{\mathrm{B}} \|u(t)\|_{L^{2}(J;\mathbb{U})} (s_2-s_1)^{\frac{1}{2}} + M_{\mathrm{B}} \|u(t)\|_{L^{2}(J;\mathbb{U})} \int_{0}^{s_1} \left\|\mathrm{U}(s_2,s)-\mathrm{U}(s_1,s)\right\|_{\mathcal{L}(\mathbb{H})} ds\nonumber\\
		&\quad + M \left(C_{f}+q^{*}C_{\xi}\right) (s_2-s_1) + \left(C_{f}+q^{*}C_{\xi}\right) \int_{0}^{s_1} \left\|\mathrm{U}(s_2,s)-\mathrm{U}(s_1,s)\right\|_{\mathcal{L}(\mathbb{H})} ds.
	\end{align}
	The right hand side of the inequality (\eqref{eq3.6}) converges to zero uniformly for $x\in \mathcal{B}_{r_0}$ as $|s_2-s_1|\rightarrow0$, since the operator $\mathrm{U}(t,s)$ is continuous in operator topology for $t\geq0$. For $t_k <t \leq t_{k+1}$, $ k\ge 1$, we can show the equicontinuity of $F_{2}$ for any $x\in \mathcal{B}_r$in the same way as above. Therefore, the image of $\mathcal{B}_{r_{0}}$ under $F_{2}$ is equicontinuous. 
	This suggests that $F_{2}\left(B_{r_{0}}\right)$ is equicontinuous.As a result, by applying the extended version of the Arzelà-Ascoli theorem, we conclude that, $F_{2}\left(B_{r_{0}}\right)$ is relatively compact set. Hence, by Lemma\ref{fixed point thm}, the operator $F_{1}+F_{2}$ possesses at least one fixed point $x \in \mathcal{B}_{r_{0}}$, which coincides with the mild solution of system \eqref{1}.\\
	\begin{rem}
		We can also show the uniqueness of the mild solution by using the contraction mapping principle with the constant $k=\max \left\{k_{1}, k_{2}\right\}<1$, where  $k_{1}$ and $k_{2}$ are defined as
		\begin{align*}
			k_{1}= M b \left(L_{f}+q^{*}L_{\xi}\right), \quad k_{2}=  \left(M^{2} N b +  M b \right)\left(L_{f}+q^{*} L_{\xi}\right) .
		\end{align*}.
	\end{rem}
	
	Our next target is to prove the approximate controllability of semilinear system\eqref{1}.
	
	\begin{thm}
		Let the assumptions (R1)-(R4), (A1)-(A3) and the conditions of theorem \ref{Existence} are true. Then, the system \eqref{1} is approximately controllable.
	\end{thm}
	Proof: From theorem\ref{Existence}, we know that for every $\lambda >0$ and $h\in \mathbb{H}$, there exists a mild solution $x_{\lambda}\in \mathcal{PC}\left([0,b],\mathbb{H}\right)$ such that
	\begin{equation}\label{6}
		x_{\lambda}(t)=	\begin{cases}
			\mathrm{U}(t,0) x(0)+\int_{0}^{t} \mathrm{U}(t,s) \big [\mathrm{B} u(s)+f(s,x_{\lambda}(s))+\int_{0}^{s} q(s-\tau) \xi(\tau, x_{\lambda}(\tau))d\tau\big ] d s, \quad 0 \leq t \leq t_{1} \\
			\mathrm{U}\left(t,t_{k}\right) x\left(t_{k}^{+}\right)+\int_{t_{k}}^{t} \mathrm{U}(t,s) \big[\mathrm{B} u(s)+f(s,x_{\lambda}(s))+\int_{0}^{s} q(s-\tau) \xi(\tau, x_{\lambda}(\tau))d\tau\big] d s, \\ \qquad\qquad t_{k}<t \leq t_{k+1},  k=1, \ldots, m,
		\end{cases}
	\end{equation}
	where
	\begin{align*}
		\begin{split}
			x\left(t_{k}^{+}\right)  =&\prod_{j=k}^{1}\left(\mathrm{I}+\mathrm{D}_{j}\right) \mathrm{U}\left(t_{j},t_{j-1}\right) x_{0} +\sum_{i=1}^{k} \prod_{j=k}^{i+1}\left(\mathrm{I}+\mathrm{D}_{j}\right) \mathrm{U}\left(t_{j},t_{j-1}\right)\left(\mathrm{I}+\mathrm{D}_{i}\right) \\
			&\qquad \times\int_{t_{i-1}}^{t_{i}} \mathrm{T}\left(t_{i}-s\right) \big[\mathrm{B}u(s)+ f(s,x_{\lambda}(s))+\int_{0}^{s} q(s-\tau) \xi(\tau, x_{\lambda}(\tau))d\tau\big] d s \\
			&+\sum_{i=2}^{k} \prod_{j=k}^{i}\left(\mathrm{I}+\mathrm{D}_{j}\right)\mathrm{T}\left(t_{j}-t_{j-1}\right) \mathrm{E}_{i-1} v_{i-1}+\mathrm{E}_{k} v_{k}. 
		\end{split}
	\end{align*}
	The control $u(s)$ is defined as
	
	\begin{align}\label{7}
		\begin{split}
			u(s)=\bigg(\sum_{k=1}^{m}\mathrm{B}^{*}\mathrm{U}^{*}(t_{k},s)\prod_{i=k+1}^{m} \mathrm{U}(t_{i},t_{i-1})^{*}\mathrm{U}(b,t_{m})^{*}\chi(t_{k-1},t_{k})\\
			+\mathrm{B}^{*}\mathrm{U}(b,s)^{*}\chi(t_{m},b)\bigg)\widehat{\varphi}_{\lambda},\\
			v_{m}=\mathrm{E}_{m}^{*}\mathrm{U}(b,t_{m})^{*}\widehat{\varphi}_{\lambda},\quad
			v_{k}=\mathrm{E}_{k}^{*}\prod_{i=k}^{m} \mathrm{U}(t_{i},t_{i-1})^{*}(\mathrm{I}+\mathrm{D}_{i}^{*})\mathrm{U}(b,t_{m})^{*}\widehat{\varphi}_{\lambda},
		\end{split}
	\end{align}
	with 
	\begin{align*}
		\widehat{\varphi}_{\lambda}=&\bigg(\lambda \mathcal{I}+\Theta_{0}^{t_{m}}+\Gamma_{t_{m}}^{b}+\widetilde{\Theta}_{0}^{t_{m}}+\widetilde{\Gamma}_{t_{m}}^{b}\bigg)^{-1}g(x_{\lambda}(.)),
	\end{align*} 
	and
	\begin{align*}
		g(x_{\lambda}(.))=&\bigg(h-\mathrm{U}\left(b,t_{m}\right) \prod_{j=m}^{1}(\mathrm{I}+\mathrm{D}_{j}) \mathrm{U}\left(t_{j},t_{j-1}\right) x_{0}\\&  \qquad-\int_{t_{m}}^{b} \mathrm{U}(b,s)\big[f(s,x_{\lambda}(s))+\int_{0}^{s} q(s-\tau) \xi(\tau, x_{\lambda}(\tau))d\tau\big] d s\\
		&\qquad-\mathrm{U}\left(b,t_{m}\right)\sum_{i=1}^{m} \prod_{j=m}^{i+1}\left(\mathrm{I}+\mathrm{D}_{j}\right) \mathrm{U}\left(t_{j},t_{j-1}\right)\left(\mathrm{I}+\mathrm{D}_{i}\right)\int_{t_{i-1}}^{t_{i}} \mathrm{U}\left(t_{i},s\right)\\ \qquad \qquad \qquad & \qquad \qquad \qquad \qquad \big[ f(s,x_{\lambda}(s))+\int_{0}^{s} q(s-\tau) \xi(\tau, x_{\lambda}(\tau))d\tau\big] d s\bigg).
	\end{align*}
	Using \eqref{6} and \eqref{7} we can easily obtain that
	\begin{align*}
		x_{\lambda}(b)-h=&\lambda \widehat{\varphi}_{\lambda}
		=\lambda \bigg(\lambda \mathcal{I}+\Theta_{0}^{t_{m}}+\Gamma_{t_{m}}^{b}+\widetilde{\Theta}_{0}^{t_{m}}+\widetilde{\Gamma}_{t_{m}}^{b}\bigg)^{-1}g(x_{\lambda}(.)).
	\end{align*}
	
	Now, by using assumptions (A2) , we find 
	\begin{align*}
		\int_{0}^{b}\left\|f(s,x_{\lambda}(s))\right\|_{\mathbb{H}}^{2} d s\leq C_{f}^{2}b, \text{ and } 
	\end{align*}
	and the boundedness of the sequence $\left\{f(.,x_{\lambda}(.)): \lambda>0\right\}$ in $\L^{2}\left([0,b];\mathbb{H}\right)$. Then there is a subsequence still denoted by $\left\{f\left(.,x_{\lambda}(.)\right)\right\}$ that weakly converges to, say $f(.)$ in $\L^{2}\left([0,b];\mathbb{H}\right)$. Similarly by using (A3), we obtain the weak convergence of  $\left\{\xi\left(.,x_{\lambda}(.)\right)\right\}$ that weakly converges to, say $\xi(.)$ in $\L^{2}\left([0,b];\mathbb{H}\right)$.  Then by Corollary 3.3 (chapter 3) \cite{li2012optimal}, we obtain
	\begin{align}\label{8}
		\left\|g\left(x_{\lambda}(.)\right)-\omega\right\|&\leq \bigg\|\int_{t_{m}}^{b} \mathrm{U}(b,s)\big[\left(f(s,x_{\lambda}(s))-f(s)\right)+\int_{0}^{s} q(s-\tau) \left(\xi(\tau, x_{\lambda}(\tau))-\xi(\tau)\right)d\tau\big] d s \nonumber\\
		&\qquad-\mathrm{U}\left(b,t_{m}\right)\sum_{i=1}^{m} \prod_{j=m}^{i+1}\left(\mathrm{I}+\mathrm{D}_{j}\right) \mathrm{U}\left(t_{j},t_{j-1}\right)\left(\mathrm{I}+\mathrm{D}_{i}\right)\int_{t_{i-1}}^{t_{i}} \mathrm{U}\left(t_{i},s\right) \nonumber\\ \qquad \qquad \qquad & \qquad \qquad \qquad \qquad \big[ \left(f(s,x_{\lambda}(s))-f(s)\right)+\int_{0}^{s} q(s-\tau) \left(\xi(\tau, x_{\lambda}(\tau))-\xi(\tau)\right)d\tau\big] d s\bigg\|\nonumber\\&\rightarrow 0,
	\end{align}
	where
	\begin{align*}
		\omega&= h-\mathrm{U}\left(b,t_{m}\right) \prod_{j=m}^{1}(\mathrm{I}+\mathrm{D}_{j}) \mathrm{U}\left(t_{j},t_{j-1}\right) x_{0}-\int_{t_{m}}^{b} \mathrm{U}(b,s)\left[f(s)-\int_{0}^{s}q(s-\tau)\xi(\tau)d\tau\right] d s\\
		&\qquad-\mathrm{U}\left(b,t_{m}\right)\sum_{i=1}^{m} \prod_{j=m}^{i+1}\left(\mathrm{I}+\mathrm{D}_{j}\right) \mathrm{U}\left(t_{j},t_{j-1}\right)\left(\mathrm{I}+\mathrm{D}_{i}\right)\int_{t_{i-1}}^{t_{i}} \mathrm{U}\left(t_{i},s\right)\left[f(s)-\int_{0}^{s}q(s-\tau)\xi(\tau)d\tau\right]  d s,
	\end{align*}
	as $\lambda \rightarrow 0^{+}$. The first term in the right hand side of the expression \ref{8} goes to zero because of the compactness of the operator $(\mathrm{Q} f)(.)=\int_{0}^{b} \mathrm{U}(.,s)f(s) d s:\L^{2}\left([0,b];\mathbb{H}\right)\rightarrow \mathcal{PC}\left([0,b],\mathbb{H}\right) $( see Lemma 4.1 and theorem 4.2 in \cite{ravikumar2020approximate})
	and the second term tends to zero by using the compactness of the operator $\mathrm{U}(t,s)$, for $t\geq0$.
	Finally we compute $\left\|	x_{\lambda}(b)-h\right\|_{\mathbb{H}}$ as
	\begin{align*}
		\left\|x_{\lambda}(b)-h\right\|=&\left\|\lambda \bigg(\lambda \mathcal{I}+\Theta_{0}^{t_{m}}+\Gamma_{t_{m}}^{b}+\widetilde{\Theta}_{0}^{t_{m}}+\widetilde{\Gamma}_{t_{m}}^{b}\bigg)^{-1}g(x_{\lambda}(.))\right\|,\\
		&\leq \left\|\lambda \bigg(\lambda \mathcal{I}+\Theta_{0}^{t_{m}}+\Gamma_{t_{m}}^{b}+\widetilde{\Theta}_{0}^{t_{m}}+\widetilde{\Gamma}_{t_{m}}^{b}\bigg)^{-1}\omega\right\|\\&\qquad+\left\|\lambda \bigg(\lambda \mathcal{I}+\Theta_{0}^{t_{m}}+\Gamma_{t_{m}}^{b}+\widetilde{\Theta}_{0}^{t_{m}}+\widetilde{\Gamma}_{t_{m}}^{b}\bigg)^{-1}\left(g\left(x_{\lambda}(.)\right)-\omega\right)\right\|.
	\end{align*}
	By estimate \eqref{8} and assumption (A1), we obtain 
	\begin{align*}
		\left\|	x_{\lambda}(b)-h\right\|_{\mathbb{H}}\rightarrow 0 \text{ as } \lambda\rightarrow 0^{+}.
	\end{align*}
	which guarantee that the system \eqref{1} is approximately controllable in $\mathbb{H}$.
	
	\section{\textbf{Application}}
	We consider the following impulsive semilinear functional heat problem on $\mathbb{H}=\mathbb{U}=L^{2}\left([0,\pi];\mathbb{R}\right)$:
	
	\begin{equation}\label{example}
		\left\{
		\begin{aligned}
			\frac{\partial}{\partial t}z(t,\zeta)&= a(t)\frac{\partial^{2}}{\partial\zeta^{2}}z(t,\zeta)+\mu(t,\zeta)+ \frac{e^{-t}z\left(t,\zeta\right)}{(9+e^{t})(1+z(t,\zeta))} +\int_{0}^{t} e^{t-s} \frac{e^{s} z\left(s,\zeta\right)}{5+z\left(s,\zeta\right)} ds,\quad \zeta\in[0,\pi],\\
			&\qquad t\in [0,1], t\ne \left\{\frac{1}{2}\right\},\\
			z(t,0)&=0=z(t,\pi), \quad t\in [0,1],\\
			\Delta z\left(\frac{1}{2},\zeta\right)&= \mathrm{D}_{1}z\left(\frac{1}{2},\zeta\right)+\mathrm{E}v_{1},\\
			\Delta z\left(1,\zeta\right)&= \mathrm{D}_{2}z\left(1,\zeta\right)+\mathrm{E}_{2}v_{2},\\
			z\left(0,\zeta\right)&= \phi(\zeta).
		\end{aligned}
		\right.
	\end{equation}
where $a:[0,1]\mapsto\mathbb{R}^{+}$, is Holder continuous function of order $0<\mathscr{k}\leq1$, that is there exists a positive constant $C_{a}$ such that 
\begin{align*}
	|a(t)-a(s)|\leq C_{a}|t-s|^{\mathscr{k}}, \text{for all} t,s\in [0,1].
\end{align*}
For $\mathbb{H}=L^{2}\left([0,\pi];\mathbb{R}\right)$, the operator $\mathrm{A}(t) g(\zeta)= a(t)g^{\prime \prime}(\zeta)$, with the domain $\mathcal{D}(\mathrm{A}(t))=\mathcal{D}(\mathrm{A})= \mathrm{H}^{2}\left([0,\pi];\mathbb{R}\right)\cap \mathrm{W}_{0}^{1,2}([0,\pi];\mathbb{R})$. We define the operator $\mathrm{A}(t)$ as $\mathrm{A}g(\zeta)=g^{\prime \prime}$, $\zeta \in [0,\pi]$, with the domain $\mathcal{D}(\mathrm{A})$. Moreover, for $t\in [0,1]$ and $g\in \mathcal{D}(\mathrm{A})$, the operator $\mathrm{A}(t)$ can be expressed as 
\begin{align*}
	\mathrm{A}(t)g=\sum_{n=1}^{\infty}(-n^{2}a(t))\langle g,w_{n} \rangle w_{n},\quad g\in \mathcal{D}(\mathrm{A}) ,\text{ for } \quad \langle g, w_n \rangle = \int_{0}^{\pi} g(\zeta) w_n 9\zeta d\zeta,
\end{align*}
where, $-n^{2} (n\in 
\mathbb{N})$ and $w_n(\zeta)= \sqrt{\frac{2}{\pi}}\sin (n\zeta)$, are the eigenvalues and the corresponding normalizes eigenfunctions of the operator $\mathrm{A}$ respectively. The operator $\mathrm{A}(t)$ satisfies all the conditions (R1)-(R4) of assumption \ref*{2.1}(see application section of \cite{ravikumar2020approximate}). Then by applying Lemma\ref{lem 2.1}, we obtain the existence of a unique evolution system $\{\mathrm{U}(t,s): 0\leq s \leq t \leq 1\}$. From Lemma\ref{lem2.2}, we observe that the evolution system $\{\mathrm{U}(t,s): 0\leq s \leq t \leq 1\}$ is compact for $t-s>0$. The evolution system $\mathrm{U}(t,s)$ can be explicitly as
\begin{align*}
	\mathrm{U}(t,s)g= \sum_{n=1}^{\infty} e^{-n^2 \int_{s}^{t}a(\tau) d\tau}\langle g, w_n \rangle w_n, \quad \text{ for each} \quad g\in \mathbb{H}. 
\end{align*}
We also have 
\begin{align*}
	\mathrm{U}(t,s)^{*}g^{*}= \sum_{n=1}^{\infty} e^{-n^2 \int_{s}^{t}a(\tau) d\tau}\langle g^{*}, w_n \rangle w_n, \quad \text{ for each} \quad g^{*}\in \mathbb{H}. 
\end{align*}

	Next, we define operator $\mathrm{B}:L^{2}\left([0,\pi];\mathbb{R}\right)\rightarrow\mathbb{H}$ such that
	\begin{align*}
		\mathrm{B}\left(u(t)\right)\left(\zeta\right)=u(t)\left(\zeta\right)=\mu\left(t,\zeta\right), \quad t\in [0,1],\quad \zeta\in [0,\pi].
	\end{align*}
	We can see, the operator $\mathrm{B}$ defined as above is a linear bounded operator. We also define $\mathrm{D}_{k}=\mathrm{E}_{k}=\mathcal{I}, \text{ for } k=1,2$.
	
	Let the function $x:J\rightarrow\mathbb{H}$ be given by
	\begin{align*}
		x(t)\left(\zeta\right)=z\left(t,\zeta\right),\quad \zeta\in [0,\pi].
	\end{align*}
  The nonlinear functions $f, \xi:[0,1]\times D\rightarrow\mathbb{H}$ is defined as 
	\begin{align*}
		f\left(t,x(t)\right)\left(\zeta\right)= \frac{e^{-t}z\left(t,\zeta\right)}{(9+e^{t})(1+z(t,\zeta))} \text{ and } \xi\left(t,x(t)\right)\left(\zeta\right)= \frac{e^{t} z\left(t,\zeta\right)}{5+z\left(t,\zeta\right)},\quad \zeta\in[0,\pi].
	\end{align*}
We can check that for $f$ and $\xi$, assumptions (A2) and (A3) are satisfied with $L_{f}=\frac{1}{10}$ , $L_{\xi}=\frac{e}{25}$, $C_{f}=\frac{1}{10}$, $ C_{\xi}= \frac{e}{5}$. We take $v_{1}=\sin(\pi t)$, $v_{2}= \cos(\pi t)$ and $q^{*}= e-1$.
	By the above settings we can transform system\eqref{example} in the abstract form as system \eqref{1}.

Since all the conditions are satisfied therefore, there exists a mild solution the system \eqref{example} and is approximately controllable. 
	
	\section{\textbf{Conclusion}}
	In this study, we explored the solution and controllability for a class of nonautonomous impulsive integro differential systems within a Hilbert space. We first established the existence of mild solutions for the system by utilizing Krasnoselskii's fixed point theorem.
 Furthermore, we demonstrated the system’s approximate controllability. To substantiate the theoretical findings, we also provided a comprehensive example. This research advances the understanding of control methods for impulsive nonlinear systems and can be extended for second order.

	\bibliographystyle{plainnat}
	\bibliography{references}
	
\end{document}